\documentclass[11pt,a4paper]{article} 
\usepackage{graphicx,amsmath,bm}
\usepackage{amssymb,amscd}
\usepackage{caption2,color}

\usepackage[dvipdfm,
           pdfstartview=FitH,
           CJKbookmarks=true,
           bookmarksnumbered=true,
           bookmarksopen=true,
          colorlinks=true,
          colorlinks=cyan,
           pdfborder=001,
           citecolor=blue%
           ]{hyperref}

\makeatletter

\@addtoreset{equation}{section} \makeatother

\newtheorem{lemma}{Lemma}[section]

\newtheorem{theorem}{Theorem}[section]
\newtheorem{corollary}{Corollary}[section]
\newcommand{\dx}{\,\mathrm{d}x}
\newcommand{\dt}{\,\mathrm{d}t}
\newcommand{\ds}{\,\mathrm{d}s}
\newcommand{\n}{\nabla}
\newcommand{\p}{\partial}
\newcommand{\norm}[1]{\lVert#1\rVert}
\newcommand{\seminorm}[1]{\lvert#1\rvert}
\newcommand{\rn}{\ensuremath{\mathbb{R}^N}}
\DeclareMathAlphabet{\mathsfsl}{OT1}{cmss}{m}{sl}

\renewcommand{\vec}[1]{\mbox{\boldmath$#1$}}
\newcommand{\oo}{\ensuremath{\Omega}}
\newcommand{\G}{\Gamma}

\newcommand{\diff}{\,\mathrm{d}}

\newcommand{\mdiv}{\,\mathrm{div}\,}
\newcommand{\mcurl}{\,\mathrm{curl}\,}

\newcommand{\md}{\mathrm{D}}

\newcommand{\vv}{\vec{V}}
\newcommand{\vw}{\vec{w}}
\newcommand{\vvv}{\vec{v}}

\newcommand{\vy}{\vec{y}}
\newcommand{\vyd}{\vec{y}_d}

\newcommand{\vu}{\vec{u}}
\newcommand{\vf}{\vec{f}}
\newcommand{\vg}{\vec{g}}
\newcommand{\vh}{\vec{h}}
\newcommand{\vphi}{\vec{\varphi}}

\newcommand{\vpsi}{\vec{\psi}}

\newcommand{\vn}{\vec{n}}
\begin{document}

\title{Optimal Shape Design for the Time-dependent\\ Navier--Stokes Flow\footnote{This work was
supported by the National Natural Science Fund of China under grant
numbers 10371096, 10671153 for ZM Gao and YC Ma.}}

\author{Zhiming Gao\thanks{Corresponding author.  School of Science, Xi'an Jiaotong University, P.O.Box 1844, Xi'an, Shaanxi, P.R.China, 710049. E--mail:dtgaozm@gmail.com.}\qquad
 Yichen Ma\footnote{School of Science, Xi'an Jiaotong University, Shaanxi, P.R.China, 710049. E-mail:\,ycma@mail.xjtu.edu.cn.}
 \qquad Hongwei Zhuang\thanks{Engineering College of Armed Police Force, Shaanxi,\,P.R.China, 710086.}}
\date{}
 \maketitle
\noindent{{\textbf{Abstract.\;}} This paper is concerned with the
problem of shape optimization of two-dimensional flows governed by
the time-dependent Navier-Stokes equations. We derive the structures
of shape gradients with respect to the shape of the variable domain
for time--dependent cost functionals by using the state derivative
with respect to the shape of the fluid domain and its associated
adjoint state. Finally we apply a gradient type algorithm to our
problem and numerical examples show that our theory is useful for
practical purpose and the proposed algorithm is feasible in low Reynolds number flow.   \\[8pt]
{{\textbf{Keywords.\;}}
shape optimization; shape derivative; gradient algorithm; material derivative; time--dependent Navier-Stokes equations.\\[8pt]
{{\textbf{AMS(2000) subject classifications.\;}}35B37, 35Q30, 49K40.

\section{Introduction}
The problem of finding the optimal design of a system governed by
the incompressible Navier-Stokes equations arises in many design
problems in aerospace, automotive, hydraulic, ocean, structural, and
wind engineering. Example applications include aerodynamic design of
automotive vehicles, trains, low speed aircraft, sails, and
hydrodynamic design of ship hulls, turbomachinery, and offshore
structures. In many cases, the flow equations do not admit
steady-state solutions, and the optimization model must incorporate
the time-dependent form of the Navier-Stokes equations.

Optimal shape design has received considerable attention already.
Early works concerning on existence of solutions and
differentiability of the quantity (such as, state, cost functional,
etc.)
 with respect to shape deformation occupied most of the 1980s
 (see \cite{ce81,delfour,piron74,Piron,zolesio}), the stabilization of structures
 using boundary variation technique has been fully addressed in \cite{delfour,Piron,zolesio}.
However, a few studies have considered the shape optimization of
time-dependent flows (see \cite{dziri02,bei97,piron01,yagi05}). Our
concern in this article is on shape sensitivity analysis of
time-dependent Navier-Stokes flow with small regularity data, and on
deriving an efficient numerical approach for the solution of
two-dimensional realizations of such problems.

In \cite{gao-robin}, we use the state derivative approach to solve a
shape optimization problem governed by a Robin problem, and in
\cite{gao06a,gao06b}, we derive the expression of shape gradients
for Stokes and Navier--Stokes optimization problem by this approach,
respectively. In this paper, we use this approach and weak implicit
function theorem to derive the structures of shape gradients with
respect to the shape of the variable domain for some given cost
functionals in shape optimization problems for time--dependent
Navier--Stokes flow with small regularity data.

This paper is organized as follows. In section 2, we briefly recall
the velocity method which is used for the characterization of the
deformation of the shape of the domain and give the definitions of
Eulerian derivative and shape derivative. We also give the
description of the shape optimization problem for the
time--dependent Navier--Stokes flow.

In section 3, we employ the weak implicit function theorem to prove
the existence of the weak Piola material derivative, and then give
the description of the shape derivative. After that, we express the
shape gradients of some typical cost functionals by introducing the
corresponding linear adjoint state systems.

 Finally in section 4, we
propose a gradient type algorithm with some numerical examples to
prove that our theory could be very useful for the practical purpose
and the proposed algorithm is efficient in low Reynolds number flow.

\section{Preliminaries and statement of the problem}\label{sec2}
\subsection{Elements of the velocity method and notations}
Domains $\oo$ don't belong to a vector space and this requires the
development of {shape calculus} to make sense of a ``derivative" or
a ``gradient". To realize it, there are about three types of
techniques: J.Hadamard \cite{ha07}'s normal variation method, the
{perturbation of the identity} method by J.Simon \cite{si80} and the
{velocity method} (see J.Cea\cite{ce81} and
J.-P.Zolesio\cite{delfour,zo79}). We will use the velocity method
which contains the others. In that purpose, we choose an open set
$D$ in $\rn$ with the boundary $\p D$ piecewise $C^k$, and a
velocity space $\vec V\in \mathrm{E}^k :=\{\vec V\in
C([0,\varepsilon];\mathcal{D}^k(\bar{D},\rn)): \vec V\cdot\vn_{\p
D}=0\;\mbox{on }\p D\}$, where $\varepsilon$ is a small positive
real number and $\mathcal{D}^k(\bar{D},\rn)$ denotes the space of
all $k-$times continuous differentiable functions with compact
support contained in $\rn$ . The velocity field
$$\vec V(s)(x)=\vec V(s,x), \qquad x\in D,\quad s\geq 0$$
belongs to $\mathcal{D}^k(\bar{D},\rn)$ for each $s$. It can
generate transformations
 $$T_s(\vec V)X=x(s,X),\quad s\geq 0,\quad X\in D$$
through the following dynamical system
\begin{equation}\label{dynamical}
  \left\{%
  \begin{array}{ll}
  \frac{\diff x}{\diff s}(s,X)=\vec V(s,x(s))\\[3pt]
  x(0,X)=X
  \end{array}%
  \right.
\end{equation}
with the initial value $X$ given. We denote the "transformed domain"
$T_s(\vec V)(\oo)$ by $\oo_s(\vec V)$ at $s\geq 0$, and also set
$\p\oo_s:=T_s(\p\oo)$.

There exists an interval $I=[0,\delta)$, $0<\delta\leq\varepsilon,$
and a one-to-one map $T_s$ from $\bar{D}$ onto $\bar{D}$ such that
\begin{itemize}
    \item [(i)] $T_0=\mathrm{I};$
    \item [(ii)] $(s,x)\mapsto T_s(x)$ belongs to $C^1(I;C^k(D;D))$ with $T_s(\p D)=\p D$;
    \item[(iii)]$(s,x)\mapsto T_s^{-1}(x)$ belongs to $C(I;C^k(D;D))$.
\end{itemize}
Such transformation are well studied in \cite{delfour}.

 Furthermore, for sufficiently small $s>0,$ the Jacobian $J_s$ is
 strictly positive:
 \begin{equation}\label{jacobian}
   J_s(x):=\det\seminorm{\md T_s(x)}=\det\md T_s(x)>0,
 \end{equation}
where $\md T_s(x)$ denotes the Jacobian matrix of the transformation
$T_s$ evaluated at a point $x\in D$ associated with the velocity
field $\vec V$. We will also use the following notation: $\md
T_s^{-1}(x)$ is the inverse of the matrix $\md T_s(x)$ , ${}^*\md
T_s^{-1}(x)$ is the transpose of the matrix $\md T_s^{-1}(x)$. These
quantities also satisfy the following lemmas.
\begin{lemma}[\cite{zolesio}]\label{lem:a}
    For any $\vec V\in E^k$, $\md T_s$ and $J_s$ are invertible. Moreover, $\md T_s$,
     $\md T_s^{-1}$ are in $C^1([0,\varepsilon];C^{k-1}(\bar{D};\mathbb{R}^{N\times N}))$, and $J_s$, $J_s^{-1}$ are in $C^1([0,\varepsilon];C^{k-1}(\bar{D};\mathbb{R}))$
\end{lemma}
\begin{lemma}[\cite{zolesio}]
  \label{lem:b} $\vphi$ is assumed to be a vector function in
  $C^1(D)^N$.
  \begin{itemize}
    \item [(1)] $\md(T^{-1}_s)\circ T_s=\md T_s^{-1}$;
    \item [(2)] $\md(\vphi\circ T^{-1}_s)=(\md\vphi\cdot\md
    T_s^{-1})\circ T^{-1}_s$;
    \item [(3)] $(\md\vphi)\circ T_s=\md(\vphi\circ T_s)\cdot\md
    T_s^{-1}.$
  \end{itemize}
\end{lemma}
%
%

Now let $J(\oo)$ be a real valued functional associated with any
regular domain $\oo$, we say that this functional has a {\bf
Eulerian derivative} at
$\oo$ in the direction $\vec V$ if the limit\\[6pt]
\begin{equation*}
\lim_{s\searrow 0}\frac{J(\oo_s)-J(\oo)}{s}:=\diff J(\oo;\vec V)
\end{equation*}
exists.

 Furthermore, if the map
 $$\vec V\mapsto\diff
J(\oo;\vec V):\;\mathrm{E}^k\rightarrow\mathbb{R}$$ is linear and
continuous, we say that $J$ is {\bf shape differentiable} at $\oo$.
In the distributional sense we have
\begin{equation}\label{pri:shaped}
    \diff J(\oo;\vec V)=\langle \n J,\vec V\rangle_{\mathcal{D}^k(\bar{D},\rn)'\times \mathcal{D}^k(\bar{D},\rn)}.
\end{equation}
 When $J$ has a Eulerian derivative, we say that $\n J$ is the {\bf shape gradient} of $J$
at $\oo$.

Before closing this subsection, we introduce the following
functional spaces which will be used throughout this paper:
\begin{eqnarray*}
  H(\mdiv,\oo):=\{\vu\in L^2(\oo)^N: \;\mdiv\vu=0\mbox{ in
  }\oo,\;\vu\cdot\vn=0\mbox{ on }\p\oo\},\\
  H^1_0(\mdiv,\oo):=\{\vu\in H^1(\oo)^N:\;\mdiv\vu=0\mbox{ in
  }\oo,\;\vu|_{\p\oo}=0\}.
\end{eqnarray*}
Given $T>0$, we introduce the notation $L^p(0,T;X)$ which denotes
the space of $L^p$ integrable functions $f$ from $[0,T]$ into the
Banach space $X$ with the norm
\begin{equation*}
\norm{f}_{L^p(0,T;X)}=\left(\int^T_0\norm{f}^p_X\diff
t\right)^{1/p},\quad 1\leq p< +\infty.
\end{equation*}
We also denote by $L^\infty(0,T;X)$ the space of essentially bounded
functions $f$ from $[0,T]$ into $X$, and is equipped with the Banach
norm
\begin{equation*}
  {\mathrm{ess}}\sup\limits_{t\in [0,T]}\norm{f(t)}_X.
\end{equation*}
\subsection{Statement of the shape optimization problem}
In two dimensions, we consider a typical problem in which a solid
body $S$ with the boundary $\p S$ is located in an external flow.
Since the flow is in an unbounded domain, we reduce the problem to a
bounded domain $D$ by introducing an artificial boundary $\p D$ on
which we set the speed flow $\vy=\vy_\infty$. $\oo:=D\backslash S$
is the effective domain with its boundary $\p\oo=\p S\cup\p D$. The
state equations of the flow can be written by the Navier--Stokes
equations in the non-dimensional form,
\begin{equation}\label{tns:nonhomo}
  \left\{%
  \begin{array}{ll}
\p_t\vy-\alpha\Delta\vy+\md\vy\cdot\vy+\n p=\vf&\quad\mbox{in}\;Q:=\oo\times (0,T),\\
\mdiv \vy=0&\quad\mbox{in}\;Q,\\
\vy=\vy_\infty&\quad\mbox{on}\;\p D\times (0,T),\\
\vy=0&\quad\mbox{on}\;\p S\times (0,T),\\
\vy(0)=\vy_0&\quad\mbox{in }\oo,\\
\int_\oo p\dx=0,&\quad\mbox{on }(0,T),\\
 \int_{\p D}
        \vy_\infty\cdot\vn\ds=0,&\quad\mbox{on }(0,T),
  \end{array}%
  \right.
\end{equation}
where the last relation is needed in view of the incompressibility
constraint $\mdiv\vy=0$, $\alpha$ stands for the inverse of the
Reynolds number whenever the variables are appropriately
nondimensionalized, $\vy, p,$ and $\vf$ are the velocity, pressure,
and the given body force per unit mass, respectively.

Our goal is to optimize the shape of the boundary $\p S$ which
minimizes a given cost functional $J$ depending on the fluid state.
The cost functional may represent a given objective related to
specific characteristic features of the fluid flow (e.g., the
deviation with respect to a given target velocity, the drag, the
vorticity, ...).

Hence, we are interested in solving the following minimization
problem
\begin{equation}\label{ns:cost}
 \min_{\oo\in\mathcal{O}} J_1(\oo)=\frac{1}{2}\int^T_0\int_{\oo}\seminorm{\vy-\vyd}^2\dx\dt,
\end{equation}
or
\begin{equation}\label{ns:cost2}
    \min_{\oo\in\mathcal{O}}
    J_2(\oo)=\frac{\alpha}{2}\int^T_0\int_{\oo}\seminorm{\mcurl\vy}^2\dx\dt,
\end{equation}
where $\vy$ is satisfied by the full Navier--Stokes system
\eqref{tns:nonhomo} and $\vy_d$ is the target velocity given by the
engineers. We also notice that the boundary $\p D$ is fixed in our
optimization problems and an example of the admissible set
${\mathcal{O}}$ is:
$$\mathcal{O}:=\left\{\oo\subset\rn:\; \p D \mbox{ is fixed},\;\int_\oo\dx=\mbox{constant}\right\}.$$

 In order to deal with the
nonhomogeneous Dirichlet boundary condition on $\p D$, let the
vectorial function $\vec h$ be the solution of
\begin{equation}
  \left\{
  \begin{array}{ll}
  \mdiv\vec h=0\quad& \mbox{in }\oo\\
  \vec h=\vy_\infty\quad&\mbox{ on  }\partial D\\
  \vec h=0\quad&\mbox{ on }\partial S,
  \end{array}
  \right.
\end{equation}
then we can choose an extension $\vec h$ with $\vec h=0$ in the body
$S$.

Now we may look for a solution of the nonhomogeneous Navier--Stokes
equations in the form
\begin{equation}\label{homo:a}
  \vy=\vec h+\tilde\vy,
\end{equation}
with $\tilde\vy$ vanishing on the boundary of the domain $\oo$.
Substituting \eqref{homo:a} in the system \eqref{tns:nonhomo}, we
find the following equations for $\tilde\vy$:
\begin{equation}\label{tns:homo}
  \left\{%
  \begin{array}{ll}
\p_t\tilde\vy-\alpha\Delta\tilde\vy+\md\tilde\vy\cdot\tilde\vy+\md\tilde\vy\cdot\vec
h
+\md\vec h\cdot\tilde\vy+\n p=\vec F&\quad\mbox{in}\;Q,\\
\mdiv \tilde\vy=0&\quad\mbox{in}\;Q,\\
\tilde\vy=0&\quad\mbox{on}\;\p \oo\times (0,T),\\
\tilde\vy(0)=\tilde\vy_0&\quad\mbox{in }\oo
  \end{array}%
  \right.
\end{equation}
where $\vec F:=\vf+\alpha\Delta\vec h-\md\vec h\cdot\vec h$ and
$\tilde\vy_0:=\vy_0-\vec h$.

 For the existence and uniqueness of the
solution of the full Navier--Stokes system (\ref{tns:homo}), we have
the following results (see \cite{temam01}).
\begin{theorem}
        \label{thm:ns}
        The domain $\oo$ is supposed to be piecewise $C^1$. We
        assume that
        \begin{eqnarray}
       \label{thm1:a} &   \vf, \p_t\vf\in L^2(0,T;H(\mdiv,D)), \\
        \label{thm1:b}& \vy_0\in H^2(D)^N\cap H^1_0(\mdiv,D),\\
       \label{thm1:c} &\vy_\infty\in H^{3/2}(\p D)^N,\qquad
        \end{eqnarray}
 and the solution of \eqref{tns:nonhomo} is unique and satisfies
$$\tilde\vy, \p_t\tilde\vy\in L^2(0,T;H^1_0(\mdiv,\oo))\cap
 L^\infty(0,T;H(\mdiv,\oo)).$$ Moreover, if $\oo$ is of class
$C^2$ and $\vf\in L^\infty(0,T;H(\mdiv,D))$, then the function
$\tilde\vy\in L^\infty(0,T;H^2(\oo)^N)$.
\end{theorem}

\section{State derivative approach}
In this section, we shall prove the main theorem using an approach
based on the differentiability of the solution of the Navier--Stokes
system \eqref{tns:homo} with respect to the variable domain. To
begin with, we use the Piola transformation to bypass the divergence
free condition and then derive a weak material derivative by the
weak implicit function theorem. Then we will derive the structure of
the shape gradients of the cost functionals by introducing the
associated adjoint state equations.
\subsection{Piola material derivative}
From now on, we assume that $\oo$ is of class $C^1$ and
\eqref{thm1:a}--\eqref{thm1:c} hold. Then we say that the function
$\tilde\vy\in L^2(0,T;H^1_0(\mdiv,\oo))$ is called a weak solution
of problem (\ref{tns:homo}) if it satisfies
\begin{equation}
    \label{tns:weak}
    \langle e(\tilde\vy),\vw\rangle=\vec 0,\qquad \vw\in
    L^2(0,T;H^1_0(\mdiv,\oo)),
\end{equation}
with $e(\tilde\vy):=(e_1(\tilde\vy),e_2(\tilde\vy)),$ $\vec
0:=(0,0)$, and
\begin{equation}
    \label{ns:a}
    \langle e_1(\tilde\vy),\vw\rangle:=\int^T_0
    \int_\oo(\p_t\tilde\vy\cdot\vw+
    \alpha\md\tilde\vy:\md\vw+\md\tilde\vy\cdot\tilde\vy\cdot\vw
    +\md\tilde\vy\cdot\vh\cdot\vw+\md\vh\cdot\tilde\vy\cdot\vw-\vec
    F\cdot\vw)\dx\dt,
\end{equation}
\begin{equation}\label{ns:a0}
\langle
e_2(\tilde\vy),\vw\rangle:=\int_\oo(\tilde\vy(0)-\tilde\vy_0)\cdot\vw(0)\dx.\hspace{7.25cm}
\end{equation}
 It must be considered that the divergence free condition is variant with respect to the use of the transformation $T_s$
  during the derivation of the shape gradient for the cost functional. Therefore, we need to introduce the
  well known Piola transformation which preserves the divergence
  free condition.
\begin{lemma}[\cite{bios}]
    The Piola transform
    \begin{eqnarray*}
        \Psi_s:\quad H(\mdiv,\oo)&\mapsto& H(\mdiv,\oo_s)\\
                             \vphi&\mapsto& ((J_s)^{-1}\md T_s\cdot\vphi)\circ T_s^{-1}
    \end{eqnarray*}
is an isomorphism.
\end{lemma}

Now by the transformation $T_s$, we consider the solution
$\tilde\vy_s$ defined on $\oo_s\times (0,T)$ of the perturbed weak
formulation
\begin{eqnarray}
  &&\nonumber\int^T_0\int_{\oo_s}(\p_t\tilde\vy_s\cdot\vw_s+\alpha\md\tilde\vy_s:\md\vw_s+\md\tilde\vy_s\cdot\tilde\vy_s\cdot\vw_s
  +\md\tilde\vy_s\cdot\vec
  h\cdot\vw_s\\\label{piola:oot2}
  &&{\hspace{6.4cm}}+\md\vh\cdot\tilde\vy_s\cdot\vw_s-\vec
  F\cdot\vw_s)\dx\dt=0,\\\label{piola:oot20}
  &&\int_{\oo_s}(\tilde\vy_s(0)-\tilde\vy_0)\cdot \vw_s(0)\dx=0,
\end{eqnarray}
for all $\vw_s\in
  L^2(0,T;H^1_0(\mdiv,\oo_s))$, and introduce $\tilde\vy^s=\Psi_s^{-1}(\tilde\vy_s),
\vw^s=\Psi_s^{-1}(\vw_s)$ defined on $Q$. Then we replace
$\tilde\vy_s$, $\vw_s$ by $\Psi_s(\tilde\vy^s), \Psi_s(\vw^s)$ in
the weak system \eqref{piola:oot2}\eqref{piola:oot20}
\begin{multline}\label{piola:oot3}
  \int^T_0\int_{\oo_s}\left[(\p_t\Psi_s(\tilde\vy^s)\cdot\Psi_s(\vw^s)+
\alpha\md(\Psi_s(\tilde\vy^s)):\md(\Psi_s(\vw^s))\right.\\
    \left. +\md(\Psi_s(\tilde\vy^s))\cdot\Psi_s(\tilde\vy^s)\cdot\Psi_s(\vw^s)
    +\md(\Psi_s(\tilde\vy^s))\cdot\vh\cdot\Psi_s(\vw^s)\right.\\
    \left.
   +\md\vh\cdot\Psi_s(\tilde\vy^s)\cdot\Psi_s(\vw^s)-\vec F\cdot\Psi_s(\vw^s)\right]\dx\dt=0,
\end{multline}
\begin{equation}
  \int_{\oo_s}[\Psi_s(\tilde\vy^s(0))-\tilde\vy_0]\cdot\Psi_s(\vw^s(0))\dx=0\hspace{6cm}
\end{equation} for all $\vw^s\in L^2(0,T;H^1_0(\mdiv, \oo))$.

Using a back transport into $\oo$ and employing Lemma \ref{lem:b},
we obtain the following weak formulation
\begin{equation}\label{es:a}
    \langle e(s,\tilde\vy^s),\vw^s\rangle=\vec 0,\qquad \forall\vw^s\in L^2(0,T;
    H^1_0(\mdiv,\oo))
\end{equation}
with the notations $e:=(e_1,e_2)$, where
\begin{multline}
  \langle e_1(s,\vec v),\vw\rangle:=\int^T_0\int_\oo
  \p_t(B(s)\vvv)\cdot(\md
  T_s\vw)\dx\dt\\
  +\alpha\int^T_0\int_\oo \md(B(s)\vvv):[\md(B(s)\vw)\cdot
  A(s)]\dx\dt
  +\int^T_0\int_\oo \md(B(s)\vvv)\cdot\vvv\cdot(B(s)\vw)\dx\dt\\
  +\int^T_0\int_\oo \md(B(s)\vvv)\cdot\vh\cdot(B(s)\vw)\dx\dt
  +\int^T_0\int_\oo \md(B(s)\vh)\cdot\vvv\cdot(B(s)\vw)\dx\dt\\
  -\int^T_0\int_\oo (\vec F\circ T_s)\cdot(\md T_s\cdot\vw)\dx\dt,
\end{multline}and
\begin{equation}
\langle e_2(s,\vec
v),\vw\rangle:=\int_\oo(B(s)\vvv(0)-\tilde\vy_0\circ T_s)\cdot(\md
T_s\vw(0))\dx\hspace{3.5cm} \end{equation}
 and
\begin{equation*}
    A(s):=J_s\md T_s^{-1}{ }^*\md T_s^{-1};\qquad B(s)\vec \tau:=J_s^{-1}\md T_s\cdot\vec \tau.
\end{equation*}
Now we are interested in the differentiability of the mapping
\begin{equation*}
s \mapsto \tilde\vy^s=\Psi_s^{-1}(\tilde\vy_s):\; [0,\varepsilon]\mapsto L^2(0,T;H^1_0(\mdiv,\oo))\\
\end{equation*}
where $\varepsilon>0$ is sufficiently small and $\tilde\vy^s$ is the
solution of the weak formulation
\begin{equation}
    \langle e(s,\vvv),\vw\rangle=\vec 0,\qquad \forall \vw\in L^2(0,T;H^1_0(\mdiv,\oo)).
\end{equation}
In order to prove the differentiability of $\tilde\vy^s$ with
respect to $s$ in a neighborhood of $s=0$, there maybe two
approaches:
\begin{itemize}
    \item [(i)] analysis of the differential quotient: $\lim\limits_{s\rightarrow 0}(\tilde\vy^s-\tilde\vy)/s$;
    \item [(ii)] derivation of the local differentiability of the solution $\tilde\vy$ associated
    to the implicit equation (\ref{tns:weak}).
\end{itemize}
We use the second approach. Since $\vf\in L^2(0,T;H(\mdiv,D))$, we
deduce that $(\vf\circ T_s-\vf)/s$ weakly converges to
$\md\vf\cdot\vec V$ in $L^2(0,T;H^{-1}(D)^N)$ as $s$ goes to zero.
Thus we can not use the classical implicit function theorem, since
it requires strong differentiability results in $H^{-1}$. Hence we
introduce the following weak implicit function theorem.
\begin{theorem}[\cite{zo79}]\label{theorem:weak}
    Let $X$, $Y'$ be two Banach spaces, $I$ an open bounded set in $\mathbb{R}$, and consider the map
    \begin{equation*}
      (s,x)\mapsto e(s,x):\;  I\times X\mapsto Y'
    \end{equation*}
    If the following hypothesis hold:
    \begin{itemize}
        \item [(i)]$s\mapsto \langle e(s,x),y\rangle$ is continuously differentiable for any $y\in Y$ and
        $(s,x)\mapsto\langle \p_s e(s,x),y\rangle$ is continuous;
        \item [(ii)] there exists $u\in X$ such that $u\in C^{0,1}(I;X)$ and $e(s,u(s))=0$, $\forall s\in I$;
        \item [(iii)] $x\mapsto e(s,x)$ is differentiable and $(s,x)\mapsto \p_x e(s,x)$ is continuous;
        \item [(iv)] there exists $s_0\in I$ such that $\p_x e(s,x)|_{(s_0,x(s_0))}$ is an isomorphism from $X$ to $Y'$,
    \end{itemize}
   the mapping
    \begin{equation*}
        s\mapsto u(s):\; I \mapsto X
    \end{equation*}
    is differentiable at $s=s_0$ for the weak topology in $X$ and its weak derivative $\dot{u}(s)$ is the solution of
    \begin{equation*}
        \langle \p_x e(s_0,u(s_0))\cdot\dot{u}(s_0),y\rangle+\langle \p_s e(s_0,u(s_0)),y\rangle=0,\quad \forall y\in Y.
    \end{equation*}
\end{theorem}
We may now state the main theorem of this section.
\begin{theorem}\label{thm:piola}
   We assume that the domain $\oo$ is piecewise $C^1$ and \eqref{thm1:a}--\eqref{thm1:c} hold,
   $\tilde\vy\in L^2(0,T;H^1_0(\mdiv,\oo))$ is the solution of the weak formulation
\eqref{tns:weak}. Then the weak {\bf Piola material derivative}
    $\dot{\tilde\vy}^P:=\p_s(\tilde\vy^s)|_{s=0}$
    exists and is characterized by the following weak formulation:
    \begin{equation}
        \langle \p_{\vec v} e(0,\vec v)|_{\vec v=\tilde\vy}\cdot{\dot{\tilde\vy}}^P,\vw\rangle
        +\langle \p_s e(0,\tilde\vy),\vw\rangle=\vec 0,\quad\forall \vw\in L^2(0,T;H^1_0(\mdiv,\oo)),
    \end{equation}
i.e.,
\begin{multline}\label{tpiola:weak}
    \int^T_0\int_\oo[\p_t\dot{\tilde\vy}^P\cdot\vw+\alpha\md\dot{\tilde\vy}^P:\md\vw+
    \md\dot{\tilde\vy}^P\cdot\tilde\vy\cdot\vw
    +\md{\tilde\vy}\cdot\dot{\tilde\vy}^P\cdot\vw+\md\dot{\tilde\vy}^P\cdot\vh\cdot\vw
    +\md\vh\cdot\dot{\tilde\vy}^P\cdot\vw]\dx\dt\\=
-\int^T_0\int_\oo[\p_t((\md\vv-\mdiv\vv)\tilde\vy)\cdot\vw+\p_t\tilde\vy\cdot\md\vv\cdot\vw]\dx\dt\\
 -\alpha\int^T_0\int_\oo\md((\md\vv-\mdiv\vv)\tilde\vy):\md\vw\dx\dt\\
 -\alpha\int^T_0\int_\oo\md\tilde\vy:[\md((\md\vv-\mdiv\vv)\vw)+\md\vw\cdot(\mdiv\vv-\md\vv-{}^*\md\vv)]\dx\dt\\
 -\int^T_0\int_\oo[\md((\md\vv-\mdiv\vv)(\tilde\vy+\vh))\cdot\tilde\vy\cdot\vw+\md(\tilde\vy+\vh)\cdot\tilde\vy
 \cdot((\md\vv-\mdiv\vv)\vw)]\dx\dt\\
-\int^T_0\int_\oo[\md((\md\vv-\mdiv\vv)\tilde\vy)\cdot\vh\cdot\vw
 -\md\tilde\vy\cdot\vh\cdot((\md\vv-\mdiv\vv)\vw)]\dx\dt\\
 +\int^T_0\int_\oo({}^*\md\vv\cdot(\vf+\alpha\Delta\vec h-\md\vec h\cdot\vec h)
 +\md(\vf+\alpha\Delta\vec h-\md\vec h\cdot\vec
 h)\cdot\vv)\cdot\vw\dx\dt\\
 +\int^T_0\int_\oo(\vf+\alpha\Delta\vec h-\md\vec h\cdot\vec
 h)\cdot(\md\vec V\cdot\vw)\dx\dt,
\end{multline}
and
\begin{equation}
\label{tpiola:weak2}
\int_\oo\dot{\tilde\vy}^P(0)\cdot\vw(0)\dx=-\int_\oo[
 (\md\vv+{}^*\md\vv-\mdiv\vv\,\mathrm{I})\cdot\tilde\vy(0)-(\md\tilde\vy_0\vec
 V+{}^*\md\vv\tilde\vy_0)]\cdot\vw(0)\dx.
\end{equation}
\end{theorem}

\noindent{\bf Proof.}\; In order to apply Theorem
\ref{theorem:weak}, we need to verify the four hypothesis of Theorem
\ref{theorem:weak} for the mapping
\begin{equation*}
(s,\vvv)\mapsto e(s,\vvv):\;[0,\varepsilon]\times
L^2(0,T;H^1_0(\mdiv,\oo))\mapsto L^2(0,T;H^1_0(\mdiv,\oo)').
\end{equation*}
To begin with, since $\oo$ is of piecewise $C^1$, the mapping
$T_s\in C^1([0,\varepsilon];C^1({D},{D}))$. Then by Lemma
\ref{lem:a}, the mapping
\begin{equation*}
    s\mapsto  \langle e_i(s,\vvv),\vw\rangle:  [0,\varepsilon]\mapsto
    \mathbb{R}\qquad (i=1,2)
\end{equation*}
is $C^1$ for any $\vvv, \vw\in L^2(0,T;H^1_0(\mdiv,\oo))$. On the
other hand, since $\vf\in L^2(0,T;$ $H(\mdiv,D))$, the mapping
$s\mapsto \vec f\circ T_s$ is only weakly differentiable in
$H^{-1}$, thus the mapping $s\mapsto e_1(s,\vvv)$ is weakly
differentiable, and then $s\mapsto e(s,\vvv)$ is weakly
differentiable.

Since we have the following identities by simple calculation,
\begin{eqnarray}
    \label{tderi:a}\frac{\diff}{\diff s}\md T_s&=&(\md\vec V(s)\circ T_s)\md T_s;\\
    \label{tderi:b}\frac{\diff}{\diff s}J_s&=&(\mdiv \vec V(s))\circ T_s\, J_s;\\
    \label{tderi:c}\frac{\diff}{\diff s}(\vf\circ T_s)&=&(\md\vf\cdot\vec V(s))\circ T_s,
\end{eqnarray}
the weak derivative of $e_i(s,\vvv)\; (i=1,2)$ can be expressed as
\begin{multline}
  \langle \p_s e_1(s,\vvv),\vw\rangle=\int^T_0\int_\oo
  [\p_t(B'(s)\vvv)\cdot(\md
  T_s\vw)+\p_t(B(s)\vvv)\cdot(\md\vv(s)\circ T_s)\cdot\md
  T_s\cdot\vw]\dx\dt\\
  +\alpha\int^T_0\int_\oo\md(B'(s)\vvv):[\md(B(s)\vw)\cdot
  A(s)]\dx\dt\\
  +\alpha\int^T_0\int_\oo \md(B(s)\vvv):[\md(B'(s)\vw)\cdot
  A(s)+\md(B(s)\vw)\cdot A'(s)]\dx\dt\\
    +\int^T_0\int_\oo[\md(B'(s)\vvv)\cdot\vvv\cdot(B(s)\vw)+\md(B(s)\vvv)\cdot\vvv\cdot(B'(s)\vw)]\dx\dt\\
    +\int^T_0\int_\oo[\md(B'(s)\vvv)\cdot\vh\cdot(B(s)\vw)+\md(B(s)\vvv)\cdot\vh\cdot(B'(s)\vw)]\dx\dt\\
        +\int^T_0\int_\oo[\md(B'(s)\vh)\cdot\vvv\cdot(B(s)\vw)+\md(B(s)\vh)\cdot\vvv\cdot(B'(s)\vw)]\dx\dt\\
        -\int^T_0\int_\oo[{ }^*\md\vec V(s)\cdot\vec F+\md\vec F\cdot\vec V(s)]\circ T_s\cdot(\md T_s\,\vw)\dx\dt\\
        -\int^T_0\int_\oo(\vec F\circ T_s)\cdot[(\md\vv(s)\circ T_s)\cdot\md
 T_s\cdot\vw]\dx\dt,
\end{multline}
and
\begin{multline}
 \langle \p_s e_2(s,\vvv),\vw\rangle=\int_\oo \{[B'(s)\vvv(0)-(\md\tilde\vy_0\cdot\vv(s))\circ
 T_s]\cdot(\md T_s\vw(0))\\
 +(B(s)\vvv(0)-\tilde\vy_0\circ T_s)\cdot[(\md\vv(s)\circ T_s)\cdot\md
 T_s\cdot\vw(0)]\}\dx,
\end{multline}
where
\begin{eqnarray*}
    &B'(s)\vec\tau:=\frac{\p}{\p s}\left[B(s)\vec\tau\right]=[\md\vec V(s)\circ T_s-(\mdiv \vec V(s)\circ T_s)\mathrm{I}]B(s)\vec \tau;\\[6pt]
    &A'(s):=\frac{\p}{\p s}A(s)=[\mdiv \vec V(s)\circ T_s-\md T_s^{-1}\md\vec V(s)\circ T_s]\,A(s)-{}^*[\md T_s^{-1}\md\vec V(s)\circ T_s\, A(s)].
\end{eqnarray*}
Obviously, the mapping $(s,\vvv)\mapsto \p_s e(s,\vvv)$ is
continuous, and when we take $s=0$, we have
\begin{eqnarray*}
  B'(0)\vec\tau=(\md\vv-\mdiv\vv\mathrm{I})\cdot\vec\tau;\\
  A'(0)=\mdiv\vv\mathrm{I}-\md\vv-{}^*\md\vv,
\end{eqnarray*}
and then
\begin{multline}\label{es:b}
 \langle \p_s
 e_1(0,\vvv),\vw\rangle=\int^T_0\int_\oo[\p_t((\md\vv-\mdiv\vv)\vvv)\cdot\vw+\p_t\vvv\cdot\md\vv\cdot\vw]\dx\dt\\
 +\alpha\int^T_0\int_\oo\md((\md\vv-\mdiv\vv)\vvv):\md\vw\dx\dt\\
 +\alpha\int^T_0\int_\oo\md\vvv:[\md((\md\vv-\mdiv\vv)\vw)+\md\vw\cdot(\mdiv\vv-\md\vv-{}^*\md\vv)]\dx\dt\\
 +\int^T_0\int_\oo[\md((\md\vv-\mdiv\vv)(\vvv+\vh))\cdot\vvv\cdot\vw+\md(\vvv+\vh)\cdot\vvv\cdot((\md\vv-\mdiv\vv)\vw)]\dx\dt\\
 +\int^T_0\int_\oo[\md((\md\vv-\mdiv\vv)\vvv)\cdot\vh\cdot\vw+\md\vvv\cdot\vh\cdot((\md\vv-\mdiv\vv)\vw)]\dx\dt\\
 -\int^T_0\int_\oo\vec F\cdot(\md\vec V\cdot\vw)\dx\dt-\int^T_0\int_\oo({}^*\md\vv\cdot\vec F+\md\vec
 F\cdot\vv)\cdot\vw\dx\dt.
\end{multline}
\begin{equation}\label{es:b2}
 \langle \p_s
 e_2(0,\vvv),\vw\rangle=\int_\oo[
 (\md\vv+{}^*\md\vv-\mdiv\vv)\cdot\vvv(0)\cdot\vw(0)-(\md\tilde\vy_0\vec
 V+{}^*\md\vv\tilde\vy_0)\cdot\vw(0)]\dx.
\end{equation} To verify (ii), we follow the same steps described in
R.Dziri\cite{dziri95} to find that the mapping $s\mapsto
\tilde\vy_s\circ T_s$ is Lipschitz continuous which is the direct
consequence of the uniqueness of the solution of the Navier--Stokes
system, i.e., Theorem \ref{thm:ns}.

It is easy to check that the mappings \begin{eqnarray*} \vvv\mapsto
e_1(s,\vvv)&:& \quad L^2(0,T;H^1_0(\mdiv,\oo))\rightarrow
L^2(0,T;H^1_0(\mdiv,\oo)')\\
\vvv\mapsto e_2(s,\vvv)&:&\quad H^1_0(\mdiv,\oo)\rightarrow
H^1_0(\mdiv,\oo)
\end{eqnarray*}
are differentiable, and the derivatives of $e_i(s,\vvv)$ with
respect to $\vvv$ in the direction $\delta\vvv$ are
\begin{multline}
  \langle\p_v e_1(s,\vvv)\cdot\delta\vvv,\vw\rangle=\int^T_0\int_\oo
  \p_t (B(s)\delta\vvv)\cdot(\md T_s\vw)\dx\dt\\
  +\alpha\int^T_0\int_\oo\md(B(s)\delta\vvv):[\md(B(s)\vw)\cdot
  A(s)]\dx\dt\\
  +\int^T_0\int_\oo[\md(B(s)\delta\vvv)\cdot\vvv\cdot(B(s)\vw)+\md(B(s)\vvv)\cdot\delta\vvv\cdot(B(s)\vw)]\dx\dt\\
  +\int^T_0\int_\oo[\md(B(s)\delta\vvv)\cdot\vh\cdot(B(s)\vw)+\md(B(s)\vh)\cdot\delta\vvv\cdot(B(s)\vw)]\dx\dt.
\end{multline}
and
\begin{equation}
\langle\p_v
e_2(s,\vvv)\cdot\delta\vvv,\vw\rangle=\int_\oo(B(s)\delta\vvv(0))\cdot(\md
T_s\cdot\vw(0))\dx.
\end{equation}
 The continuity of $(s,\vvv)\mapsto \p_v e_i(s,\vvv)$ is easy to
check. Moreover,
\begin{multline}\label{es:deltav}
   \langle\p_v
   e_1(0,\vvv)\cdot\delta\vvv,\vw\rangle=\int^T_0\int_\oo
   [\p_t(\delta\vvv)\cdot\vw+\alpha\md(\delta\vvv):\md\vw
   +\md(\delta\vvv)\cdot\vvv\cdot\vw\\\md\vvv\cdot\delta\vvv\cdot\vw
   +\md(\delta\vvv)\cdot\vh\cdot\vw+\md\vh\cdot\delta\vvv\cdot\vw]\dx\dt,
\end{multline}
\begin{equation}\label{es:deltav2}
\langle\p_v
   e_2(0,\vvv)\cdot\delta\vvv,\vw\rangle=\int_\oo
   \delta\vvv(0)\cdot\vw(0)\dx.
\end{equation} Furthermore, $\delta\vvv\rightarrow \p_v
e(0,\vvv)\cdot\delta\vvv$ is an isomorphism which follows from the
uniqueness and existence of the Navier--Stokes system, i.e., Theorem
\ref{thm:ns}. Indeed, we assume that $\tilde\vy_1, \tilde\vy_2$ are
two solutions of the Navier--Stokes system \eqref{tns:homo}, and
$\tilde\vy_i\; (i=1,2)$ satisfies the weak formulation
\eqref{tns:weak}. It is obvious that
$\hat\vy=\tilde\vy_1-\tilde\vy_2$ satisfies
\begin{equation}
  \int^T_0\int_\oo[\p_t\hat{\vy}\cdot\vw+\alpha\md\hat\vy:\md\vw+\md\hat\vy\cdot\vh\cdot\vw+\md\vh\cdot\hat\vy\cdot\vw
  +\md\hat\vy\cdot\tilde\vy_1\cdot\vw+\md\tilde\vy_2\cdot\hat\vy\cdot\vw]\dx\dt=0,
\end{equation}
and
\begin{equation}
  \int_\oo\hat\vy(0)\cdot\vw(0)\dx=0.
\end{equation} Now let $\vw=\hat\vy$, we can follow the proof of the unique solvability of the unsteady Navier--Stokes equations (see Temam \cite{temam01}) and
obtain
$$\seminorm{\hat\vy(t)}^2\leq 0,\qquad\forall t\in [0,T].$$
Thus $\tilde\vy_1=\tilde\vy_2$. Similar a priori estimates hold for
$\delta\vvv$ and the uniqueness of the solution of the system
\eqref{es:deltav}\eqref{es:deltav2} is obtained.

Finally, all the hypothesis are satisfied by (\ref{es:a}), we can
apply Theorem \ref{theorem:weak} to (\ref{es:a}) and then use
\eqref{es:b}, \eqref{es:b2}, \eqref{es:deltav} and
\eqref{es:deltav2} to obtain \eqref{tpiola:weak} and
\eqref{tpiola:weak2}.\hfill $\square$

\subsection{Shape derivative}
In this subsection, we will characterize the shape derivative
$\tilde\vy'$, i.e., the derivative of the state $\tilde\vy$ with
respect to the shape of the variable domain.
\begin{theorem}\label{tns:shaped}
   Under the assumption of Theorem \ref{thm:ns} and moreover assume that $\oo$ is of class $C^2$, $\tilde\vy\in L^\infty(0,T;H^2(\oo)^N\cap H^1_0(\mdiv,\oo))$ solves the weak formulation
    \eqref{tns:weak} and $\tilde\vy_s$ solves the perturbed weak formulation \eqref{piola:oot2}\eqref{piola:oot20} in $\oo_s\times (0,T)$, then the {\bf shape derivative}
    $$\tilde\vy':=\lim_{s\rightarrow 0}\frac{\tilde\vy_s-\tilde\vy}{s}$$ exists and is characterized as the solution of
\begin{equation}\label{tns:shape}
    \left\{
    \begin{array}{lll}
        \p_t\tilde\vy'-\alpha\Delta\tilde\vy'
        +\md\tilde\vy'\cdot\tilde\vy+\md\tilde\vy\cdot\tilde\vy'+\md\tilde\vy'\cdot\vh+\md\vh\cdot\tilde\vy'+\n p'=0
        \;&\mbox{in }Q\\
        \mdiv \tilde\vy'=0 &\mbox{in  }Q\\
        \tilde\vy'=-(\md\tilde\vy\cdot\vn)\vv_n&\mbox{on
        }\p S\times (0,T)\\
        \tilde\vy'=0&\mbox{on }\p D\times (0,T)\\
        \tilde\vy'(0)=0&\mbox{in }\oo.
    \end{array}
    \right.
\end{equation}
\end{theorem}

\noindent{\bf Proof.}\; Since $\oo$ is of class $C^2$ and $\vec
V\in\mathrm{E}^2$, $\oo_s$ has the same regularity than $\oo$ for
any $s\in (0,\epsilon)$, then $\tilde\vy_s\in
L^\infty(0,T;H^2(\oo_s)^N)$ satisfies the following weak formulation
\begin{equation}\label{tstokes:c}
    \int^T_0\int_{\oo_s}(\alpha\md\tilde\vy_s:\md\vw+\md\tilde\vy_s\cdot\tilde\vy_s\cdot\vw
    +\md\tilde\vy_s\cdot\vh\cdot\vw+\md\vh\cdot\tilde\vy_s\cdot\vw-\vec
    F\cdot\vw)\dx\dt=0,
\end{equation}
\begin{equation}\label{tstokes:d}
  \int_{\oo_s}\tilde\vy_s(0)\cdot\vw(0)\dx=0
\end{equation} for any $\vw\in L^2(0,T;H^1_0(\mdiv,\oo_s))$.
Moreover, we have $\p_t\tilde\vy_s\in L^2(0,T;H^1_0(\mdiv,\oo))$.

To begin with, we introduce the following Hadamard formula (see
\cite{delfour,zolesio})
\begin{equation}\label{hadamard}
 \frac{\diff{}}{\diff s}\int_{\oo_s}g(s,x)\dx=\int_{\oo_s}
 \frac{\p g}{\p s}(s,x)\dx+\int_{\p\oo_s} g(s,x)\,\vec
 V\cdot\vn_s\diff\G_s,
 \end{equation}
 for a sufficiently smooth functional
$g:[0,\tau]\times\rn\rightarrow\mathbb{R}$.

Now we set a function $\vphi\in\mathcal{D}(Q)^{N+1}$ and
$\mdiv\vphi(x,t)=0$ in $\oo$ for a.e. $t\in (0,T)$. Obviously when
$s$ is sufficiently small, $\vphi(t)$ belongs to the sobolev space
$H^1_0(\mdiv,\oo_s)\cap H^2(\oo_s)^N$ for a.e. $t\in (0,T)$. Hence
we can use (\ref{hadamard}) to differentiate \eqref{tstokes:c},
\eqref{tstokes:d} with $\vw=\vphi$,
\begin{multline*}
    \int^T_0\int_\oo(\p_t\tilde\vy'+\alpha\md\tilde\vy':\md\vphi+\md\tilde\vy'\cdot\tilde\vy\cdot\vphi+\md\tilde\vy\cdot\tilde\vy'\cdot\vphi
    +\md\tilde\vy'\cdot\vh\cdot\vphi+\md\vh\cdot\tilde\vy'\cdot\vphi)\dx\dt\\
    +\int^T_0\int_{\p\oo} (\alpha\md\tilde\vy:\md\vphi+\md\tilde\vy\cdot\tilde\vy\cdot\vphi+\md\tilde\vy\cdot\vh\cdot\vphi
    +\md\vh\cdot\tilde\vy\cdot\vphi-\vec F\cdot\vphi)\vv_n\diff
    s\dt=0,
\end{multline*}$$
    \int_\oo\tilde\vy'(0)\cdot\vphi(0)\dx+\int_{\p\oo}\tilde\vy_s(0)\cdot\vphi(0)\vec
    V_n\diff s=0.$$
Since $\vphi$ has a compact support, the boundary integrals vanish.
Using integration by parts, we obtain
\begin{equation}
    \int^T_0\int_\oo(\p_t\tilde\vy'-\alpha\Delta\tilde\vy'+\md\tilde\vy'\cdot\tilde\vy+\md\tilde\vy\cdot\tilde\vy'
    +\md\tilde\vy'\cdot\vh
    +\md\vh\cdot\tilde\vy')\cdot\vphi\dx\dt=0,
\end{equation}
and
\begin{equation}
 \int_\oo\tilde\vy'(0)\cdot\vphi(0)\dx=0.
\end{equation}
 Then there exists some distribution $p'$ such that
$$\p_t\tilde\vy'-\alpha\Delta\tilde\vy'+\md\tilde\vy'\cdot\tilde\vy
+\md\tilde\vy\cdot\tilde\vy'+\md\tilde\vy'\cdot\vh+\md\vh\cdot\tilde\vy'=-\n
p'$$ in the distributional sense in $Q$ and $\tilde\vy'(0)=0$ in
$\oo$ since $\vphi(0)$ is arbitrary.

Now we recall that for each sufficient small $s$,
$\Psi_s^{-1}(\tilde\vy_s)$ belongs to the Sobolev space
$H^1_0(\mdiv,\oo)$, then we can deduce that its material derivative
vanishes on the boundary $\p S$. Thus we obtain the shape derivative
of $\tilde\vy$ at the boundary $\p S$,
\begin{equation*}
    \tilde\vy'=-\md\tilde\vy\cdot\vv,\qquad\mbox{on  }\p S\times
    (0,T).
\end{equation*}
Since $\tilde\vy|_{\p S\times (0,T)}=0,$ we have $\md\tilde\vy|_{\p
S\times (0,T)}=\md\tilde\vy\cdot\vn^*\vn,$ and then
\begin{equation*}
    \tilde\vy'=-(\md\tilde\vy\cdot\vn)\vv_n\qquad\mbox{on  }\p S\times (0,T).
\end{equation*}
Since $\p D$ is fixed, we obtain $\tilde\vy'=0$ on the boundary $\p
D\times (0,T)$. \hfill$\square$

The shape derivative $\vy'$ of the solution $\vy$ of the original
Navier--Stokes system (\ref{tns:nonhomo}) is given by
$\tilde\vy'=\vy'$, then we obtain the following corollary by
substituting $\tilde\vy'=\vy'$ and $\tilde\vy=\vy-\vh$ into
(\ref{tns:shape}).
\begin{corollary}
    The shape derivative $\vy'$ of the solution $\vy$ of (\ref{tns:nonhomo}) exists and satisfies the following system
    \begin{equation}\label{tns:shapederivative}
        \left\{
        \begin{array}{ll}
            \p_t\vy'-\alpha\Delta\vy'+\md\vy'\cdot\vy+\md\vy\cdot\vy'+\n p'=0 &\quad\mbox{in }Q;\\
            \mdiv\vy'=0&\quad\text{in }Q;\\   \vy'=(-\md\vy\cdot\vn)\vv_n&\quad\mbox{on }\p S\times
            (0,T)\\
            \vy'=0&\quad\mbox{on }\p D\times
            (0,T)\\
            \vy'(0)=0&\quad\mbox{in }\oo.
        \end{array}
        \right.
    \end{equation}
\end{corollary}

\subsection{Adjoint state system and gradients of the cost functionals}
This subsection is devoted to the computation of the shape gradients
for the cost functionals $J_1(\oo)$ and $J_2(\oo)$ by the adjoint
method.

For the cost functional
$J_1(\oo)=\int^T_0\int_\oo\frac{1}{2}\seminorm{\vy-\vy_d}^2\dx\dt$,
we have
\begin{theorem}\label{thm:a}
    Let $\oo$ be of class $C^2$, $\vy_d\in L^\infty(0,T;L^2(D)^N)$, and $\vec V\in \mathrm{E}^2$, the shape gradient $\n J_1$ of the cost functional $J_1(\oo)$ can be expressed as
    \begin{equation}\label{nsa:gradient}
        \n J_1=\left[\frac{1}{2}(\vy-\vy_d)^2+\alpha(\md\vy\cdot\vn)\cdot(\md\vec v\cdot\vn)\right]\vn,
    \end{equation}
    where the adjoint state $\vec v$ satisfies the following linear adjoint system
    \begin{equation}\label{adjoint:a}
\left\{
\begin{array}{lll}
    -\p_t\vec v-\alpha\Delta\vec v-\md\vec v\cdot\vy+{ }^*\md\vy\cdot\vec v+\n q=\vy-\vy_d,&\qquad\mbox{in  }Q\\
    \mdiv\vec v=0,&\qquad\mbox{in  }Q\\
    \vec v=0,&\qquad\mbox{on  }\p\oo\times (0,T)\\
    \vec v(T)=0,&\qquad\mbox{in }\oo.
\end{array}
\right.
    \end{equation}
\end{theorem}
\noindent{\bf Proof.}\; Since $J_1(\oo)$ is differentiable with
respect to $\vy$, and the state $\vy$ is shape differentiable with
respect to $s$, i.e., the shape derivative $\vy'$ exists, we obtain
Eulerian derivative of $J_1(\oo)$ with respect to $s$,
\begin{equation}\label{a:b}
    \diff J_1(\oo;\vec V)=\int^T_0\int_\oo (\vy-\vy_d)\cdot\vy'\dx\dt
    +\int^T_0\int_{\p\oo} \frac{1}{2}\seminorm{\vy-\vy_d}^2\vec V_n\diff
    s\dt
\end{equation}
by Hadamard formula (\ref{hadamard}).

By Green formula, we have the following identity
\begin{multline}\label{a:a}
    \int^T_0\int_\oo [(\p_t\vy'-\alpha\Delta\vy'+\md\vy'\cdot\vy+\md\vy\cdot\vy'+\n p')\cdot\vw-\mdiv\vy'\pi]\dx\dt
    \\=\int^T_0\int_\oo[(-\p_t\vw-\alpha\Delta\vw-\md\vw\cdot\vy+{ }^*\md\vy\cdot\vw+\n\pi)\cdot\vy'-p'\mdiv\vw]\dx\dt\\
+\int^T_0\int_{\p\oo} (\vy'\cdot\vw)(\vy\cdot\vn)\diff
s\dt+\int^T_0\int_{\p\oo}
(\alpha\md\vw\cdot\vn-\pi\vn)\cdot\vy'\diff
s\dt\\+\int^T_0\int_{\p\oo} (p'\vn-\alpha\md\vy'\vn)\cdot\vw\diff
s\dt+\int_\oo (\vy'(T)\cdot\vw(T)-\vy'(0)\cdot\vw(0))\dx.
\end{multline}
Now we define $(\vec v,q)$ to be the solution of (\ref{adjoint:a}),
use (\ref{tns:shapederivative}) and set $(\vw,\pi)= (\vec v,q)$ in
(\ref{a:a}) to obtain
\begin{equation}
\int^T_0\int_\oo (\vy-\vy_d)\cdot\vy'\dx\dt=-\int^T_0\int_{\p
S}(\alpha\md\vec v\cdot\vn-q\vn)\cdot\vy'\diff s\dt.
\end{equation}
Since $\vy'=(-\md\vy\cdot\vn)\vec V_n$ on the boundary $\p S$ and
$\mdiv\vy'=0$ in $\oo$, we obtain the Eulerian derivative of
$J_1(\oo)$ from (\ref{a:b}),
\begin{equation}
    \diff J_1(\oo;\vec V)=\int^T_0\int_{\p S} \left[\frac{1}{2}\seminorm{\vy-\vy_d}^2+
    \alpha\left(\md\vy\cdot\vn\right)\cdot(\md\vec v\cdot\vn)\right]\vec V_n\diff s\dt.
\end{equation}
Since the mapping $\vec V\mapsto \diff J_1(\oo;\vec V)$ is linear
and continuous, we get the expression (\ref{nsa:gradient}) for the
shape gradient $\n J_1$ by (\ref{pri:shaped}).\hfill $\square$

For another typical cost functional
$J_2(\oo)=\frac{\alpha}{2}\int^T_0\int_\oo\seminorm{\mcurl\vy}^2\dx\dt$,
we have the following theorem.
\begin{theorem}
    Let $\oo$ be of class $C^2$ and $\vec V\in \mathrm{E}^2,$ the cost functional $J_2(\oo)$ possesses the shape gradient $\n J_2$ which can be expressed as
    \begin{equation}\label{nsb:gradient}
        \n J_2=\alpha\left[\frac{1}{2}\seminorm{\mathrm{curl}\,\vy}^2+(\md\vy\cdot\vn)\cdot(\md\vec v\cdot\vn-\mathrm{curl}\,\vy\wedge\vn)\right]\vn,
    \end{equation}
    where the adjoint state $\vec v$ satisfies the following linear adjoint system
    \begin{equation}\label{adjoint:b}
\left\{
\begin{array}{lll}
   -\p_t\vec v-\alpha\Delta\vec v-\md\vec v\cdot\vy+{ }^*\md\vy\cdot\vec v+\n q=-\alpha\Delta\vy,&\qquad\mbox{in  }Q\\
    \mdiv\vec v=0,&\qquad\mbox{in  }Q\\
    \vec v=0,&\qquad\mbox{on  }\p\oo\times (0,T)\\
    \vec v(T)=0,&\qquad\mbox{in }\oo.
    \end{array}
\right.
    \end{equation}
\end{theorem}
\noindent{\bf Proof.}\; The proof is similar to that of Theorem
\ref{thm:a}. Using Hadamard formula (\ref{hadamard}) for the cost
functional $J_2$, we obtain the Eulerian derivative
\begin{equation}\label{b:b}
    \diff J_2(\oo;\vec V)=\alpha\int^T_0\int_\oo\mcurl\vy\cdot\mcurl\vy'\dx\dt
    +\int^T_0\int_{\p\oo}\frac{\alpha}{2}\seminorm{\mcurl\vy}^2\vec V_n\diff s\dt.
\end{equation}
Then, we define $(\vec v,q)$ to be the solution of
(\ref{adjoint:b}), use (\ref{tns:shapederivative}) and set
$(\vw,\pi)= (\vec v,q)$ in (\ref{a:a}) to obtain
\begin{equation}\label{b:c}
\alpha\int^T_0\int_\oo\Delta\vy\cdot\vy'\dx\dt=\int^T_0\int_{\p S}
\alpha(\md\vec v\cdot\vn)\cdot\vy'\diff s\dt.
\end{equation}
Applying the following vectorial Green formula
\begin{multline*}
\int_\oo
(\vphi\cdot\Delta\vpsi+\mcurl\vphi\cdot\mcurl\vpsi+\mdiv\vphi\mdiv\vpsi)\dx\\=\int_{\p\oo}
(\vphi\cdot(\mcurl\vpsi\wedge\vn)+\vphi\cdot\vn\mdiv\vpsi)\diff s
\end{multline*}
for the vector functions $\vy$ and $\vy'$, we obtain
\begin{equation}\label{b:d}
\int^T_0\int_\oo(\mcurl\vy\cdot\mcurl\vy'+\Delta\vy\cdot\vy')\dx\dt=\int^T_0
\int_{\p S}(\mcurl\vy\wedge\vn)\cdot\vy'\diff s\dt
\end{equation}
Combining (\ref{b:b}), (\ref{b:c}) with (\ref{b:d}), we obtain the
Eulerian derivative
\begin{equation*}
\diff J_2(\oo;\vec V)=\int^T_0\int_{\p
S}\alpha\left[\frac{1}{2}\seminorm{\mathrm{curl}\,\vy}^2+(\md(\vy-\vg)\cdot\vn)\cdot(\md\vec
v\cdot\vn-\mathrm{curl}\,\vy\wedge\vn)\right]\vec V_n\diff s\dt.
\end{equation*}
Finally we arrive at the expression (\ref{nsb:gradient}) for the
shape gradient $\n J_2$.\hfill $\square$
\section{Gradient algorithm and numerical simulation}
In this section, we will give a gradient type algorithm and some
numerical examples in two dimensions to prove that our previous
methods could be very useful and efficient for the numerical
implementation of the shape optimization problems for the unsteady
Navier--Stokes flow. For the sake of simplicity, we only consider
the cost functional $J(\oo)=\int^T_0\int_\oo
\seminorm{\vy-\vy_d}^2\dx\dt.$

\subsection{A gradient type algorithm} As
we have just seen, the general form of the Eulerian derivative is
\begin{equation*}
  \diff J(\oo;\vec V)=\int^T_0\int_{\p S} \n J\cdot \vec V\ds\dt,
\end{equation*}
where $\n J$ denotes the shape gradient of the cost functional $J$.
Ignoring regularization, a descent direction is found by defining
\begin{equation}
  \vec V=-h_k\n J
\end{equation}
and then we can update the shape $\oo$ as
\begin{equation}
  \oo_k=(\mathrm{I}+h_k\vec V)\oo
\end{equation}
where $h_k$ is a descent step at $k$-th iteration.

There are also other choices for the definition of the descent
direction. Since the gradient of the functional has necessarily less
regularity than the parameter, an iterative scheme like the method
of descent deteriortates the regularity of the optimized parameter.
We need to project or smooth the variation into $H^1(\oo)^2$. Hence,
the method used in this paper is to change the scalar product with
respect to which we compute a descent direction, for instance,
$H^1(\oo)^2$. In this case, the descent direction is the unique
element $\vec d\in H^1(\oo)^2$ such that at a fixed time $t\in
[0,T]$ and for every $\vec V\in H^1(\oo)^2,$
\begin{equation}\label{reg}
  \int_\oo\md\vec d :\md \vec V\dx=\int_{\p S} \n J\cdot \vec V\ds.
\end{equation}
  The computation of $\vec d$ can also be interpreted as a regularization
  of the shape gradient, and the choice of $H^1(\oo)^2$ as space of
  variations is more dictated by technical considerations rather
  than theoretical ones.

The resulting algorithm can be summarized as follows:
\begin{itemize}
    \item [(1)] Choose an initial shape $\oo_0$, i.e., choose an initial shape of $\p S$ since $\p D$ is fixed in our problem;
    \item [(2)] Compute the state system \eqref{tns:nonhomo} and adjoint state system \eqref{adjoint:a}, then
    we can evaluate the descent direction $\vec d_k$ by using (\ref{reg})
    with $\oo=\oo_k;$
    \item[(3)] Set $\oo_{k+1}=(\mathrm{Id}-h_k\vec d_k) \,\oo_k,$ where $h_k$
    is a small positive real number.
\end{itemize}

The choice of the descent step $h_k$ is not an easy task. Too big,
the algorithm is unstable; too small, the rate of convergence is
insignificant. In order to refresh $h_k$, we compare $h_k$ with
$h_{k-1}$. If $(\vec d_k,\vec d_{k-1})_{H^1}$ is negative, we should
reduce the step; on the other hand, if $\vec d_k$ and $\vec d_{k-1}$
are very close, we increase the step. In addition, if reversed
triangles are appeared when moving the mesh, we also need to reduce
the step.

In our algorithm, we do not choose any stopping criterion. A
classical stopping criterion is to find that whether the shape
gradients $\n J$ in some suitable norm is small enough. However,
since we use the continuous shape gradients, it's hopeless for us to
expect very small gradient norm because of numerical discretization
errors. Instead, we fix the number of iterations. If it is too
small, we can restart it with the previous final shape as the
initial shape.
\subsection{Numerical examples}
 To illustrate the theory, we
want to solve the following minimization problem
\begin{equation}\label{exam:fun}
    \min_{\oo}\frac{1}{2}\int^1_0\int_{\oo}\seminorm{\vy-\vyd}^2\dx\dt
\end{equation}
subject to
\begin{equation}\label{exam:state}
  \left\{%
  \begin{array}{ll}
\p_t\vy-\alpha\Delta\vy+\md\vy\cdot\vy+\n p=\vf&\quad\mbox{in}\;\oo\times(0,1)\\
\mdiv \vy=0&\quad\mbox{in}\;\oo\times (0,1)\\
\vy=0&\quad\mbox{on}\;\p S\times (0,1)\\
\vy=\vy_\infty&\quad\mbox{on}\;\p D\times (0,1)\\
\vy(0)=0&\quad\mbox{in }\oo.
  \end{array}%
  \right.
\end{equation}
Where $D:=\{(x,y)\in\mathbb{R}^2:\;x^2+y^2\leq 0.64\},$ and the
shape of the body $S$ is to be optimized. We choose the velocity
$\vy_\infty=(0.15y,-0.15x)^T$ and the body force $\vf=(f_1,f_2)^T$:
\begin{multline*}
 f_1=-\frac{45x}{31\sqrt{x^2+y^2}}+\frac{\alpha\, t\, y (15 x^2+15 y^2-1)}{5(x^2+y^2)^{3/2}}\\
 +\frac{1}{25}t^2\, x\left(-46-25x^2-25y^2-\frac{1}{x^2+y^2}+\frac{12}{\sqrt{x^2+y^2}}+60\sqrt{x^2+y^2}\right);
\end{multline*}
\begin{multline*}
 f_2=-\frac{45y}{31\sqrt{x^2+y^2}}-\frac{\alpha\, t\,x (15 x^2+15 y^2-1)}{5(x^2+y^2)^{3/2}}\\
 +\frac{1}{25}t^2\,
 y\left(-46-25x^2-25y^2-\frac{1}{x^2+y^2}+\frac{12}{\sqrt{x^2+y^2}}+60\sqrt{x^2+y^2}\right).
\end{multline*}

The target velocity $\vyd$ is determined by the data $\vf,
\vy_\infty$ and the target shape of the domain $\oo$. Our aim is to
recover the shape of $S$ which is a circle: $\p
S=\{(x,y):\;x^2+y^2=0.04\}$.

The Navier--Stokes system \eqref{tns:nonhomo} and the adjoint system
\eqref{adjoint:a} are discretized by using a mixed finite element
method. Time discretization is effected using the backward Euler
method and we assume that the time interval $[0,1]$ is divided into
equal intervals of duration $\Delta t=0.05$. Spatial discretization
is effected using the Taylor--Hood pair \cite{hood74} of finite
element spaces on a triangular mesh, i.e., the finite element spaces
are chosen to be continuous piecewise quadratic polynomials for the
velocity and continuous piecewise linear polynomials for the
pressure. Our numerical solutions are obtained under FreeFem++
\cite{hecht} and we run the program on a home PC.

We choose the initial shape of $S$ to be elliptic:
$\{(x,y):\;{x^2}/{9}+{y^2}/{4}={1}/{25}\}$, and the initial finite
element mesh was shown in \autoref{fig0}.

\begin{figure}[!htbp]
  \centering
  \includegraphics[width=2.1in]{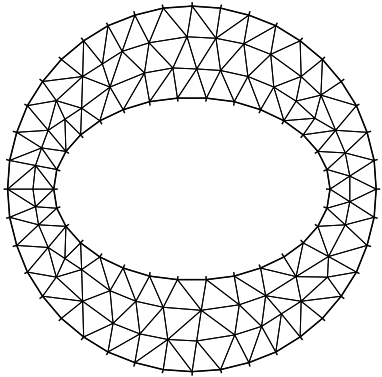}
  \caption{Initial mesh with 125 nodes.\label{fig0}}
\end{figure}

\begin{figure}[!htbp]
\centering
  \includegraphics[width=.8\textwidth]{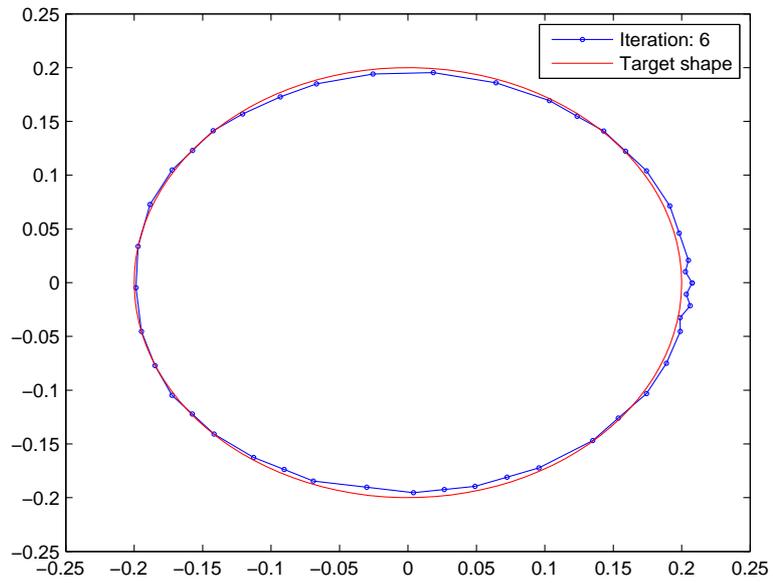}
  \caption{$\alpha=0.1$, CPU time: 124.531 s.\label{fig1}}
\end{figure}

  \autoref{fig1}---\autoref{fig3} give the comparison
between the target shape with iterated shape for the viscosity
coefficients $\alpha=0.1, 0.01$ and $0.001$, respectively. In case
of $\alpha=0.1,0.01$, we have fine results in \autoref{fig1} and
\autoref{fig2}. Unfortunately, we can not get a nice reconstruction
for $\alpha=0.001$ as in \autoref{fig3}.

\autoref{fig4} represents the fast convergence of the cost
functional for the various viscosity coefficients $\alpha=0.1, 0.01$
and $0.001$.


\begin{figure}[!htbp]
\centering
  \includegraphics[width=.8\textwidth]{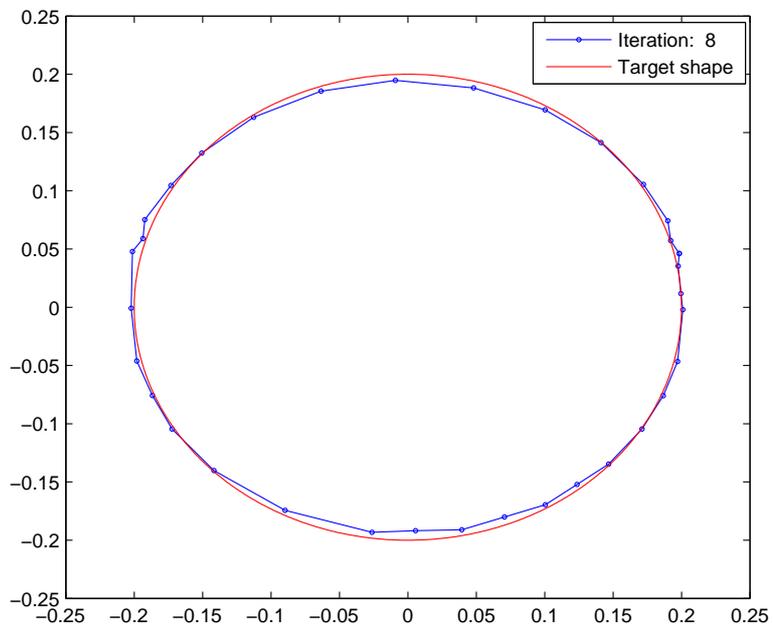}
  \caption{$\alpha=0.01$, CPU time: 120.125 s.\label{fig2}}
\end{figure}
\begin{figure}[!htbp]
\centering
  \includegraphics[width=.8\textwidth]{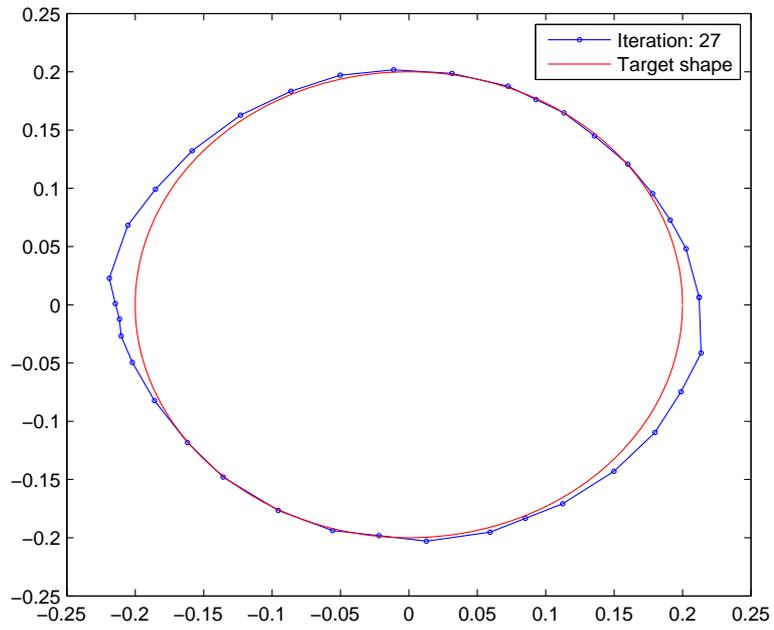}
  \caption{$\alpha=0.001$, CPU time: 622.813 s.\label{fig3}}
\end{figure}

\begin{figure}[!htbp]
\centering
  \includegraphics[width=.89\textwidth]{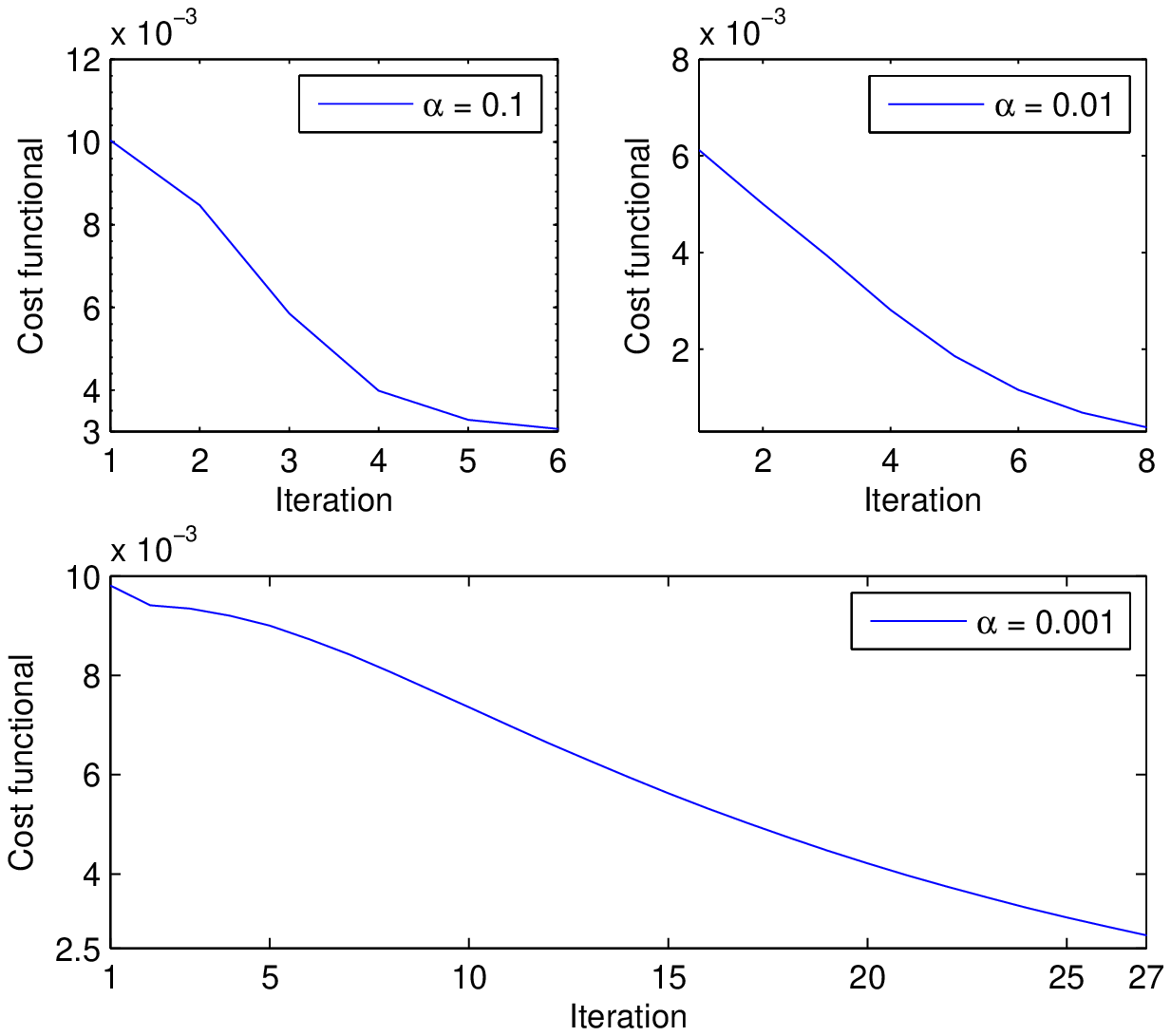}
  \caption{Convergence history for $\alpha=0.1, 0.01$ and $0.001$.\label{fig4}}
\end{figure}

\section{Conclusion}
In this paper, the shape optimization in the two dimensional
time--dependent Navier--Stokes flow has been presented. We employed
the weak implicit function theorem to obtain the existence of the
weak Piola material derivative, then we gave the description of the
shape derivative. Hence we derived the structures of shape gradients
with respect to the shape of the variable domain for some
time--dependent cost functionals by introducing the associated
adjoint state system. A gradient type algorithm is effectively used
for the minimization problem in various Reynolds number flows.
Further research is necessary on efficient implementations for very
large Reynolds numbers and real problems in the industry.

\end{document}